\documentclass[12pt]{amsart}
\usepackage{amssymb}
\usepackage{enumerate}
\usepackage{amscd}
\usepackage{graphicx}
\makeatletter
\@namedef{subjclassname@2010}{%
  \textup{2010} Mathematics Subject Classification}
\makeatother

\newtheorem{thm}{Theorem}[section]
\newtheorem{cor}[thm]{Corollary}

\newtheorem{prop}[thm]{Proposition}

\theoremstyle{definition}

\newtheorem{rem}[thm]{Remark}
\newtheorem{exa}[thm]{Example}

\numberwithin{equation}{section}

\newcommand{\pr}{\indent{\em Proof. \ }}
\newcommand{\EQ}{\begin{equation}}
\newcommand{\EN}{\end{equation}}

\newcommand{\sTr}{\mbox{\scriptsize \rm Tr}}
\newcommand{\Tr}{\mbox{\rm Tr}}

\newcommand{\F}{\mathbb{F}}

\begin{document}

\baselineskip=17pt
\title{On Kloosterman sums over finite fields of characteristic $3$
\footnote{This work has been partially supported by  the
Russian fund of fundamental researches (number of project
15 - 01 - 08051).}}

\author{L. A. Bassalygo and V.A. Zinoviev}
\address{Kharkevich Institute for Problems of Information
Transmission of the Russian Academy of Sciences,\\
Russia, 127994, Moscow, GSP-4, B. Karetnyi per. 19}
\email{bass@iitp.ru, zinov@iitp.ru}

\date{}

\begin{abstract}
We study the divisibility by $3^k$ of Kloosterman sums $K(a)$ over
finite fields of characteristic $3$.  We give a new recurrent
algorithm for finding the largest $k$, such that $3^k$ divides the
Kloosterman sum $K(a)$. This gives a new simple test for zeros
of such Kloosterman sums.
\end{abstract}



\maketitle

\section{Introduction}

Let $\F = \F_{q}$ be a field of characteristic $p$
of order $q=p^m$, where  $m \geq 2$ is an integer and let
$\F^* = \F \setminus \{0\}$. By $\F_p$ denote the field, consisting
of $p$ elements. For any element $a \in \F^*$ the {\em Klosterman sum}
can be defined as
\EQ\label{eq:1.1}
K(a)~=~\sum_{x \in \F} \omega^{\sTr(x+a/x)},
\EN
where $\omega = \exp{2\pi i/p}$ is a primitive $p$-th root of
unity and
\EQ\label{eq:1.2}
Tr(x)~=~x + x^p + x^{p^2} + \cdots + x^{p^{m-1}}.
\EN
Recall that under $x^{-i}$ we understand  $x^{p^m-1-i}$, avoiding by
this way a division into $0$.

Kloosterman sums are used for solutions of equations over finite
fields \cite{lid}, in the theory of error correcting codes \cite{lac}, for
studying and constructing of bent and hyperbent functions \cite{ch3,he2} and so
on. Surely, any Kloosterman sum $K(a)$ for a given $a$ can be
found directly computing its value for every element of the field,
but this method requires a large amount of computations, which is
a multiple of the size of the field. Hence more simple methods of
computations of Kloosterman sums are quite interesting. The values of
characteristic $p \in \{2,3\}$ are especially interesting in connection with
the number of $q$-rational points of some elliptic curves \cite{kat,lac,lis,mo1}.
Divisibility of binary Kloosterman sums by numbers $8, 16, \ldots, 256$ and
computations of such sums modulo some numbers was considered in
papers \cite{ch1,ch2,gee,go1,go2,go3,he1,lis,lim,mo2}.
Divisibility of ternary Klosterman sums $K(a)$ by $9$ and by $27$ was
considered in \cite{go2,go3,go4,gee,lis,lim}. Furthermore, in \cite{go2,go3,go4}
a complete characterization of $K(a)$ modulo $9$, $18$ and $27$ was obtained for
any $m$ and $a$.

In the recent papers \cite{ba0,ba1} we provide two algorithms to
find the maximum divisor of $K(a)$ of type $2^k$. Similar results
for the case $p=3$ have been announced in \cite{ba2} and here we
prove these results. In particular, we give a simple test of
divisibility of $K(a)$ by $27$. We suggest also a new recursive
algorithm of finding the largest divisor of $K(a)$ of the type
$3^k$ which needs at every step the limited number of arithmetic
operations in $\F$. For the case when $m = g\,h$ we derive the
exact connection between the divisibility by $3^k$ of $K(a)$ in
$\F_{3^g}$, $a \in \F_{3^g}$, and the divisibility by $3^{k'}$ of
$K(a)$ in $\F_{3^{gh}}$.

\section{Known results}

In this section we state the known results  about Kloosterman
sums $K(a)$ \cite{lis,mo1} and elliptic curves $E(a)$ \cite{eng,men}
over finite fields $\F$ of characteristic $3$. Our interest is the
divisibility of such sums by the maximal possible number of type
$3^k$ (i.e. $3^k$ divides  $K(a)$, but  $3^{k+1}$ does not divide
$K(a)$; in addition, when $K(a) = 0$ we assume that $3^m$  divides $K(a)$,
but $3^{m+1}$ does not divide; recall that $q=3^m=|\F|$).

For a given $\F$ and any $a \in \F^*$ define the elliptic curve $E(a)$
as follows:
\EQ\label{eq:1.1}
E(a)~=~\{(x,y) \in \F \times \F:~y^2 = x^3 + x^2 - a\}.
\EN
The set of $\F$-rational points of the curve $E(a)$ over $\F$ forms
a finite abelian group, which can be represented as a direct product
of a cyclic subgroup $G(a)$ of order $3^t$ and a certain subgroup
$H(a)$ of some order $s$ (which is not multiple to $3$): $E(a) = G(a) \times H(a)$,
such that
\[
|E(a)|~=~3^t \cdot s
\]
for some integers $t \geq 2$ and $s \geq 1$ (see \cite{eng}), where
$s \not \equiv 0 \pmod{3}$.

Moisio \cite{mo1} showed that
\EQ\label{eq:2.1}
|E(a)|~=~3^m + K(a),
\EN
where $|A|$ denotes the cardinality of a finite set $A$ (earlier the same
result was obtained in \cite{kat} for the curve $y^2+xy+ay=x^3$). Therefore
a Kloosterman sum $K(a)$ is divisible by $3^t$, if and only if the
number of points of the curve $E(a)$ is divisible by $3^t$. Lisonek \cite{lis}
observed, that $|E(a)|$ is divisible by $3^t$, if and only if
the group $E(a)$ contains an element of order $3^t$.

Since $|E(a)|$ is divisible by $|G(a)|$, which is equal to $3^t$, then
generator elements of $G(a)$ and only these elements are of order
$3^t$.

Let $Q = (\xi,*) \in E(a)$. Then the point $P=(x,*) \in
E(a)$, such that $Q = 3 P$ exists, if and only if the equation
\[
x^9 - \xi x^6 + a(1-\xi)x^3 -a^2(a+\xi)~=~0.
\]
has a solution in $\F$ (see \cite{eng}). This equation is equivalent to equation
\EQ\label{eq:2.2}
x^3 - \xi^{1/3} x^2 + (a(1-\xi))^{1/3}x - (a^2(a+\xi))^{1/3}~=~0.
\EN
The equation (\ref{eq:2.2}) is solvable in $\F$ if and only if (see, for example,
\cite{ah1})
\EQ\label{eq:2.3}
\Tr\left(\frac{a\sqrt{\xi^3 + \xi^2 - a}}{\xi^3}\right)~=~0~.
\EN
Since the point $(a^{1/3},\, a^{1/3})$ of $E(a)$ has the order $3$, and hence belongs
to $G(a)$, then solving the recursive equation
\EQ\label{eq:2.4}
x_i^3 - x_{i-1}^{1/3} x_i^2 + (a(1-x_{i-1}))^{1/3}x_i -
(a^2(a+x_{i-1}))^{1/3}~=~0,~~i=0,1,...
\EN
with initial value $x_0 = a^{1/3}$, we obtain
that the point $(x_i,*) \in G(a)$ for $i=0,1,\ldots, t-1$, and the
point $(x_{t-1},*)$ is a generator element of $G(a)$. Such
algorithm of finding of cardinality of $G(a)$  was
given in \cite{ah1}.

Similar method was presented in our previous papers \cite{ba0,ba1} for finite fields of
characteristic $2$. Besides, some another results have been obtained
in \cite{ba0,ba1} for the case $p=2$. Our purpose here is to generalize these results for
finite fields of characteristic $3$.

\section{New results}

We begin with a simple result. It is known \cite{ah1,gee,lim}, that $9$
divides $K(a)$ if and only if $Tr(a) = 0$. In this case $a$ can be
presented as follows: $a = z^{27} - z^9$, where $z \in \F$,
and, hence $x_0 = a^{1/3} = z^9 - z^3$ (see (\ref{eq:2.4})). We
found the expression for the next element $x_1$, namely:
\[
x_1 = (z^4-1)(z^3-1)z^2
\]
and, therefore, from the condition
(\ref{eq:2.3}), the following result holds.

\begin{prop}\label{stat:3.1}
Let $a \in \F^*$  and $Tr(a) = 0$, i.e. $a$ can be presented in
the form: $a=z^{27} - z^9$. Then $x_0 = z^9 - z^3$,~
$ x_1 = (z^4-1)(z^3-1)z^2$,
and, therefore,  $K(a)$ is divisible by $27$, if and only if
\EQ\label{eq:3.1}
\Tr\left(\frac{z^5(z-1)(z+1)^7}{(z^2+1)^3}\right)\;=\;0,
\EN
\end{prop}

This condition (\ref{eq:3.1}) is more compact than the corresponding
condition from the papers \cite{go3,go4}, where it is proven that $K(a)$
is divisible by $27$, if $\Tr(a) = 0$ and
\[
2 \sum_{1 \leq i\leq j \leq m-1} a^{3^i + 3^j}~+~\sum_{1 \leq i<
j< k \leq m-1} a^{3^i+3^j+3^k}~=~0.
\]
Emphasize once more, that similar conditions permit in \cite{go2,go3,go4}
to find all values of $K(a)$ modulo $9$, $18$ and $27$, while the
condition (\ref{eq:3.1}) gives only divisibility of $K(a)$ by $27$.

Similar to the case $p=2$ \cite{ba0,ba1}, we give now also
another algorithm to find the maximal divisor of $K(a)$
of the type $3^t$, which requires at every step  the limited
number of arithmetic operations in $\F$.

Let $a \in \F^*$ be an arbitrary element and let
$u_1, u_2, \ldots, u_{\ell}$ be a sequence of elements of $\F$,
constructed according to the following recurrent relation
(compare with (\ref{eq:2.4})):
\EQ\label{eq:3.2}
u_{i+1}~=~\frac{(u_i^3 - a)^3~+~a u_i^3}{(u_i^3 - a)^2},~~i=1,2,
\ldots,
\EN
where $(u_1,*) \in E(a)$ and
\EQ\label{eq:3.3}
\Tr\left(\frac{a \sqrt{u_1^3 + u_1^2 -a}}{u_1^3}\right)~\neq~0~.
\EN
Then the following result is valid.

\begin{thm}\label{theo:3.1}
Let $a \in \F^*$ and let $u_1, u_2, \ldots, u_{\ell}$ be a sequence of
elements of $\F$, which satisfies the recurrent relation (\ref{eq:3.2}),
where the element $u_1$ satisfies (\ref{eq:3.3}) and $(u_1,*) \in E(a)$.
Then there exists an integer $k \leq m$ such that one of the two
following cases takes place:\\
(i) either $u_k=a^{1/3}$, but the all previous elements $u_i$ are not equal to $a^{1/3}$;\\
(ii) or $u_{k+1} = u_{k+1+r}$ for a certain $r$ and the all elements $u_i$ are
different for $i < k + 1 + r$.\\
In the both cases the Kloosterman sum $K(a)$ is divisible by $3^{k}$
and is not divisible by $3^{k+1}$.
\end{thm}

\pr
Let $a \in \F^*$ and let $u_1, u_2, \ldots, u_{\ell}$ be a sequence of
elements of $\F$, which satisfies the recurrent relation (\ref{eq:3.2}),
where the element $u_1$ satisfies (\ref{eq:3.3}) and the point
$P_1=(u_1,*)$ belongs to $E(a)$.
Assume that $E(a)$ has the order $3^t \cdot s$, where $s$ is prime to $3$.
We have to show that $k=t$.

Denote $P_i=(u_i,*)$. Since $P_1=(u_1,*)$ belongs to $E(a)$, it
follows from the addition operation in the additive abelian group $E(a)$
(Table 2.3 in \cite{eng}) and from (\ref{eq:3.2}), that for $i\geq 2$ all points
$P_i$ belongs to $E(a)$ and $P_i=3^{i-1}P_1$ for $i\geq 2$.

There are only two possibilities: either $P_1 \in G(a)$, or
$P_1 \in E(a)\setminus G(a)$.

First consider the case $P_1\in G(a)$. We claim that the condition (\ref{eq:3.3})
implies that $P_1$ is a generating element of the (cyclic) group $G(a)$.
Indeed, assume that it is not the case. Then it means that there is the point
$Q\in G(a)$ such that $P_1=3\,Q$. Assuming that $Q=(x,*)$ and using the addition
operation in $E(a)$ \cite{eng}, we arrive to the following equation for $x$:
\[
x^3-u_1^{1/3}x^2+(a(1-u_1))^{1/3}x-(a^2(a+u_1))^{1/3}=0.
\]
As we already mentioned in Section 2, this equation has a solution, if and only if
\[
\Tr\left(\frac{a \sqrt{u_1^3 + u_1^2 -a}}{u_1^3}\right)~=~0,
\]
that contradicts to (\ref{eq:3.3}). We conclude that $P_1$ is a generating point
of $G(a)$ and, therefore, has the order $3^t$. This means that the point
$P_t=(u_t,*)=3^{t-1}P_1$ is of the order $3$. Since there are exactly two points
in $E(a)$ of the order $3$ \cite{eng}, namely, the points $(a^{1/3}, \pm a^{1/3})$,
it means that for any $i \leq t-1$ we have that $u_i \neq a^{1/3}$. Therefore,
$k=t$ and $K(a)$ is divisible by $3^t$.

Now consider the case when $P_1 \in E(a)\setminus G(a)$. Then the order $d$ of
the point $3^tP_1=(u_{t+1},*)$ divides $s$. The point $d\,P_1$ belongs to
the cyclic group $G(a)$, and it is a generating element, since the equality
$dP_1=3Q$ for some $Q \in E(a)$ implies the equality $P_1=3Q'$, that contradicts
to (\ref{eq:3.3}). Therefore the order of the point $P_1$ is equal to $d\cdot 3^t$
and $K(a)$ is divisible by $3^t$.

Denote by $r$ the least integer, such that
$d$ divides $3^r-1$ or $3^r+1$.  Then we have the following equalities:
\\in the first case
\[
3^{t+r} \cdot P_1~=~3^t \cdot P_1
\]
and in the second case
\[
3^{t+r} \cdot P_1~=~-\,3^t \cdot P_1.
\]
In the both cases we obtain, that $u_{t+1} = u_{t+1+r}$ and our sequence
$u_1, u_2, \ldots, u_{\ell}$ becomes
periodic with a period $r$, starting from the elemnt $u_{t+1}$.
\qed

\begin{rem}
It is clear, that, for the case (ii) of Theorem 1, this algorithm needs $k+r$
computations of values
\[
x^3-a + \frac{a\,x^3}{(x^3-a)^2}
\]
for finding the
largest divisor of $K(a)$ of the type $3^{k}$. Besides, the following
lower bound for the number of $\F$-rational points of the curve $E(a)$
is valid:
\[
|E(a)|~\geq~ 3^k(2\,r + 1)
\]
and, respectively, the following upper bound for the value of Kloosterman
sum $K(a)$ takes place:
\[
K(a) ~\leq ~ 3^m~-~3^k(2\,r+1).
\]
\end{rem}

Directly from Theorem \ref{theo:3.1} we obtain the following necessary
and sufficient condition for an element $a \in \F^*$ to be a zero of
the Kloosterman sum $K(a)$.

\begin{cor}
Let $a \in \F^*$ and $u_1, u_2, \ldots, u_{\ell}$ be a sequence of
elements of $\F$ of the order $|\F|=3^m$, which satisfies the recurrent relation
(\ref{eq:3.2}), where the element $u_1$ satisfies (\ref{eq:3.3}).
Then  $K(a) = 0$, if and only if $u_m = a^{1/3}$, and
$u_i \neq a^{1/3}$ for all $1 \leq i \leq m-1$.
\end{cor}

\begin{exa}\label{ex:3.1}
Suppose the field $\F$ of order $3^5$ is generated by $\phi(x) =
x^5 + x^4 + x^2 + 1$ and its root $\alpha$ is a primitive element
of $\F$. Take $a=\alpha^{31}$. Following to Statement \ref{stat:3.1},
present $a$ as $a=z^{27}-z^9$. We find that
$z \in \{\alpha^{16}, \alpha^{106}, \alpha^{231}\}$. The corresponding
possible values of $x_1$ are $\alpha^7, \alpha^{19}$ and
$\alpha^{105}$, respectively.  For all these values
of $z$ the condition (\ref{eq:3.1}) is satisfied. We conclude that
$K(a)$ is divisible by $27$. Choosing $x_1=\alpha^7$ and solving
the cubic equation (\ref{eq:2.4}), we obtain three solutions for
$x_2$, namely, $x_2 \in \{\alpha^{138}, \alpha^{196}, \alpha^{237}\}$.
Choose $x_2=\alpha^{138}$. Then the condition (\ref{eq:2.3})
(with $\xi=x_2$) is not valid:
\[
\Tr\left(\frac{a\sqrt{x_2^3+x_2^2-a}}{x_2^3}\right)~=~\Tr(\alpha^{202})~=~1~.
\]
It means that we find the exact divisor $3^k$ of $K(\alpha^{31})$ and
the maximal cyclic subgroup $G(\alpha^{31})$ of the curve $E(\alpha^{31})$
is of the order $27$. The $x$-th coordinates of all points $(x,*)$ of
the cyclic group $G(\alpha^{31})$ are presented below as a graph, which
gives all possible nine sequences $x_0 =a^{1/3}, x_1, x_2$.

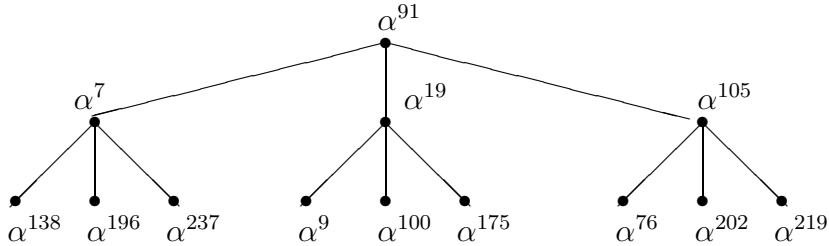
\begin{picture}(400,130)(0,0)
\put(200,90){\circle*{4}}
\put(197,95){$\alpha^{91}$}

\put(200,90){\line(-4,-1){111}}
\put(90,60){\circle*{4}}
\put(82,65){$\alpha^{7}$}

\put(200,90){\line(0,-1){30}}
\put(200,60){\circle*{4}}
\put(207,65){$\alpha^{19}$}

\put(200,90){\line(4,-1){115}}
\put(320,60){\circle*{4}}
\put(318,65){$\alpha^{105}$}

\put(90,60){\line(-1,-1){32}}
\put(60,30){\circle*{4}}
\put(57,15){$\alpha^{138}$}

\put(90,60){\line(0,-1){30}}
\put(90,30){\circle*{4}}
\put(87,15){$\alpha^{196}$}

\put(90,60){\line(1,-1){32}}
\put(120,30){\circle*{4}}
\put(117,15){$\alpha^{237}$}

\put(320,60){\line(1,-1){32}}
\put(350,30){\circle*{4}}
\put(347,15){$\alpha^{219}$}

\put(320,60){\line(0,-1){30}}
\put(320,30){\circle*{4}}
\put(317,15){$\alpha^{202}$}

\put(320,60){\line(-1,-1){32}}
\put(290,30){\circle*{4}}
\put(287,15){$\alpha^{76}$}

\put(200,60){\line(-1,-1){32}}
\put(170,30){\circle*{4}}
\put(167,15){$\alpha^{9}$}

\put(200,60){\line(0,-1){30}}
\put(200,30){\circle*{4}}
\put(197,15){$\alpha^{100}$}

\put(200,60){\line(1,-1){32}}
\put(230,30){\circle*{4}}
\put(227,15){$\alpha^{175}$}

\end{picture}

As the condition (\ref{eq:2.3}) for all elements of the last third level (numeration of
the levels $0,1, \ldots$ is starting from the top) is not satisfied, then
$k=3$, and we conclude that $K(\alpha^{31})$ is divisible by $27$
(indeed, $K(\alpha^{31}) = 27$). Note
that all points $P_i = (x_i,y_i)$, corresponding to the $i$-th level of
the graph satisfy the condition $3^{i+1} P_i = \mathcal{O}$, where
$\mathcal{O}$ is the identity element of the group $E(a)$.

Now start from the down, choosing the point $u_1=\alpha^{159}$, which
satisfies (\ref{eq:3.3}). Then using (\ref{eq:3.2}), we obtain the sequence
\[
\alpha^{159}, \alpha^{15}, \alpha^{44}, \alpha^{162}, \alpha^{162}, \ldots .
\]
We conclude that $k=3$ (see Theorem \ref{theo:3.1}) and $K(a)$ is divisible
by $3^3$ (here $r=1$). The
choice $u_1=\alpha^{193}$ results in the following sequence:
\[
\alpha^{193}, \alpha^{199}, \alpha^{50}, \alpha^{197}, \alpha^{223}, \alpha^{197}
\alpha^{223}, \ldots .
\]
We again conclude that $k=3$ (and here $r=2$).
\end{exa}

\medskip

Assume now that the field $\F_q$ of order $q = 3^m$ is embedded into
the field $\F_{q^n}$~($n \geq 2$), and $a$ is an element of $\F_q^*$.
Recall that
\[
\Tr_{q^n \rightarrow q}(x)~=~x+x^q+x^{q^2}+\ldots +
x^{q^{n-1}},~~x \in \F_{q^n}.
\]
For any elements
$a \in \F_q$ and $b \in \F_{q^n}$ define
\begin{eqnarray*}
e(a) ~=~\omega^{\sTr(a)},~~
e_n(b) ~=~\omega^{\sTr(\sTr_{q^n \rightarrow q}(b))},
\end{eqnarray*}
where $\omega$ is a primitive $3$-th root of unity.
For a given $a \in \F_q^*$ it is possible to consider the
following two Kloosterman sums:\\
\[
K(a)~=~\sum_{x \in \F_q} e\left(x + \frac{a}{x}\right),~~~
K_n(a)~=~\sum_{x \in \F_{q^n}} e_n\left(x + \frac{a}{x}\right).
\]

Denote by $H(a)$ the maximal degree of $3$, which divides
$K(a)$, and by $H_n(a)$ the maximal degree of $3$, which
divides $K_n(a)$. Recall that in the case, when $K(a) = 0$ over
$\F_q$, where $q=3^m$, we assume that $3^m$ divides $K(a)$, but
$3^{m+1}$ does not divide. There
exists a simple connection between $H(a)$ and $H_n(a)$.

\begin{thm}\label{theo:3.2}
Let $n = 3^h\cdot s$, ~$n \geq 2$,~$s \geq 1$, where $3$ and $s$
are mutually prime, and $a \in \F^*_q$. Then
\[
H_n(a)~=~H(a)~+~h.
\]
\end{thm}

The proof follows from two simple statements.

\begin{prop}\label{stat:3.1}
Let $h=0$ (that is $n$ and $3$ are coprime). Then
\[
H_n(a)~=~H(a).
\]
\end{prop}

\pr
By definition of the trace we have for any element $x \in \F_{q}$
\begin{eqnarray*}
\Tr(\Tr_{q^n \rightarrow q}(x))&=&\Tr(x+x^q+x^{q^2}+\cdots +x^{q^{n-1}})=\\
&=&\Tr(x)+\Tr(x^q)+\Tr(x^{q^2})+\cdots +\Tr(x^{q^{n-1}})=\\
&=&n\,\Tr(x)=\\
&=&\pm\,\Tr(x),
\end{eqnarray*}
where the last equality follows, since $n$ and $3$ are coprime. Therefore,
$\Tr_{q^n \rightarrow p}(x)=0$ for any $x\in \F_q$, if and only if
$\Tr_{q \rightarrow p}(x)=0$. And, since $a, a^{1/3} \in \F_q$, then
all solutions of the equation (\ref{eq:2.4}) belong to $\F_q$, that gives
the statement.
\qed

\begin{prop}\label{stat:3.2}
Let $n=3$ and $a \in \F^*_q$. Then
\[
H_3(a)=H(a) + 1.
\]
\end{prop}

\pr
It is known \cite{car} that
\EQ\label{eq:5.16}
K_3(a)~=~K(a)^3-3K(a)^2+3K(a)-3qK(a).
\EN
Assume that $K(a)$ is divisible by $3^k$. Then it is easy to see
that $K_3(a)$ is divisible by $3^{k+1}$ and is not divisible by
$3^{k+2}$.
\qed

From Theorem \ref{theo:3.2}, recalling that equality $K(a)=0$
over $\F_q$ means divisibility of $K(a)$ by $q$, we immediately
obtain the following known result \cite{lim}.

\begin{cor}
Let $a \in \F^*_q$ and $n \geq 2$. Then $K_n(a)$ is not equal to zero.
\end{cor}

\end{document}